\documentclass[graybox]{svmult}

\usepackage{type1cm}        % activate if the above 3 fonts are
\usepackage{makeidx}         % allows index generation
\usepackage{graphicx}        % standard LaTeX graphics tool
\usepackage{multicol}        % used for the two-column index
\usepackage[bottom]{footmisc}% places footnotes at page bottom

\usepackage{newtxtext}       % 
\usepackage{newtxmath}       % selects Times Roman as basic font

\makeindex             % used for the subject index
\def\cl {\nonumber \\}
\def\el {\nonumber }

\begin{document}

%\title*{A Localized Spectral-Element Reduced Basis Method for Incompressible Flow Problems}
%\title*{A Spectral Element Reduced Basis Method in Parametric CFD}
\title*{A Spectral Element Reduced Basis Method for Navier-Stokes Equations with Geometric Variations}
\titlerunning{SEM-RB for NSE with Geometric Variations}

\author{Martin W. Hess, Annalisa Quaini, and Gianluigi Rozza}
\institute{Martin W. Hess \at SISSA mathLab, International School for Advanced Studies, via Bonomea 265, I-34136 Trieste, Italy, \email{mhess@sissa.it}
\and Annalisa Quaini \at Department of Mathematics, University of Houston, Houston, Texas 77204-3008, USA \email{quaini@math.uh.edu}
\and Gianluigi Rozza \at SISSA mathLab, International School for Advanced Studies, via Bonomea 265, I-34136 Trieste, Italy, \email{gianluigi.rozza@sissa.it}}
\maketitle

\abstract*{We consider the \emph{Navier-Stokes} equations in a channel with a narrowing of varying height.
The model is discretized with high-order spectral element \emph{ansatz} functions, resulting in $6372$ degrees of freedom.
The steady-state snapshot solutions define a reduced order space through a standard POD procedure.
The reduced order space allows to accurately and efficiently evaluate the steady-state solutions for different geometries.
In particular, we detail different aspects of implementing the reduced order model in combination with a spectral element discretization.
It is shown that an expansion in element-wise local degrees of freedom can be combined with a reduced order modelling 
approach to enhance computational times in parametric many-query scenarios.}

\abstract{We consider the \emph{Navier-Stokes} equations in a channel with a narrowing of varying height.
The model is discretized with high-order spectral element \emph{ansatz} functions, resulting in $6372$ degrees of freedom.
The steady-state snapshot solutions define a reduced order space through a standard POD procedure.
The reduced order space allows to accurately and efficiently evaluate the steady-state solutions for different geometries.
In particular, we detail different aspects of implementing the reduced order model 
in combination with a spectral element discretization.
It is shown that an expansion in element-wise local degrees of freedom can be combined with a reduced order modelling 
approach to enhance computational times in parametric many-query scenarios.}

%c_f_bnd.shape  (2160, 15)
%c_f_p.shape  (1620, 15)
%c_f_int.shape  (2592, 15)
% 2160+1620+2592 = 6372

\section{Introduction and Motivation}

Spectral element methods (SEM) use high-order polynomial ansatz functions to solve partial differential equations (PDEs)
in all fields of science and engineering, see, e.g., \cite{Hess:Sherwin:2005, Hess:CHQZ1, Hess:CHQZ2, Hess:PateraSEM, Hess:Herrero2013132, Hess:Taddei}
and references therein for an overview. 
Typically, an exponential error decay under p-refinement is observed, which can provide an enhanced accuracy over standard finite element methods at the same computational cost.
In the following, we assume that the discretization error is much smaller than the model reduction error,
small enough not to interfere with our results.
 In general, this needs to be established with the use of suitable error estimation and adaptivity techniques. 

We consider the flow through a channel with a narrowing of variable height. 
A reduced order model (ROM) is computed from a few high-order SEM solves, which accurately approximates the high-order solutions for the parameter range of interest, 
i.e., the different narrowing heights under consideration. 
Since the parametric variations are affine, a mapping to a reference domain is applied without further interpolation techniques.
The focus of this work is to show how to use simulations arising from the SEM solver Nektar++ \cite{Hess:nektar} in a ROM context.
%The open-source framework ITHACA-SEM is developed to make the ROM algorithms freely available.
In particular, the multilevel static condensation of the high-order solver is not applied, but the ROM projection works with the system matrices in local coordinates. 
See \cite{Hess:Sherwin:2005} for further details. 
This is in contrast to our previous work \cite{Hess:ENUMATH17_me}, since numerical experiments 
have shown that the multilevel static condensation is inefficient in a ROM context.
Additionally, we consider affine geometry variations. With SEM as discretization method, we use global approximation functions
for the high-order as well as reduced-order methods. 
The ROM techniques described in this paper are implemented in open-source project ITHACA-SEM\footnote{\url{https://github.com/mathLab/ITHACA-SEM}}.

%Reduced order modeling replaces the high-dimensional discrete model with a low-dimensional reduced order model (ROM).
%The ROM ansatz functions are snapshots of the high-dimensional model at a few parameter values.
%Typically, an exponential decay of the approximation error between high-order model and ROM can be observed, when increasing the number of ROM ansatz functions.

The outline of the paper is as follows. In Sec.~2, the model problem is defined and the geometric variations are introduced. 
Sec.~3 provides details on the spectral element discretization, while Sec.~4 describes the model reduction approach and shows
the affine mapping to the reference domain. 
Numerical results are given in Sec.~5, while Sec.~6 summarizes the work and points out future perspectives.

\section{Problem Formulation}

Let $\Omega \in \mathbb{R}^2$ be the computational domain.
Incompressible, viscous fluid motion in spatial domain $\Omega$ over a time interval $(0, T)$
is governed by the incompressible
 \emph{Navier-Stokes} equations with vector-valued velocity $\mathbf{u}$, scalar-valued
pressure $p$, kinematic viscosity $\nu$ and a body forcing $\mathbf{f}$:
%\eqref{Hess:NSE0} - \eqref{Hess:NSE1}:
\begin{eqnarray}
\frac{\partial \mathbf{u}}{\partial t} + \mathbf{u} \cdot \nabla \mathbf{u} &=& - \nabla p + \nu \Delta \mathbf{u} + \mathbf{f}, \label{Hess:NSE0} \\
\nabla \cdot \mathbf{u} &=& 0.
\label{Hess:NSE1}
\end{eqnarray}
\noindent Boundary and initial conditions are prescribed as 
\begin{eqnarray}
\mathbf{u} &=& \mathbf{d} \quad \text{ on } \Gamma_D \times (0, T), \\
\nabla \mathbf{u} \cdot \mathbf{n} &=& \mathbf{g} \quad \text{ on } \Gamma_N \times (0, T), \\
\mathbf{u} &=& \mathbf{u}_0 \quad \text{ in } \Omega \times 0,
\label{Hess:NSE_boundaryCond}
\end{eqnarray}

\noindent with $\mathbf{d}$, $\mathbf{g}$ and $\mathbf{u}_0$ given and $\partial \Omega = \Gamma_D \cup \Gamma_N$, $\Gamma_D \cap \Gamma_N = \emptyset$.
The \emph{Reynolds} number $Re$, which characterizes the flow \cite{Hess:Holmes}, depends on $\nu$,
a characteristic velocity $U$, and a characteristic length $L$:
\begin{equation}\label{eq:re}
Re = \frac{UL}{\nu}.
\end{equation}

%\begin{eqnarray}
% Re = \frac{UL}{\nu}.
% \label{Hess:Re_visc}
%\end{eqnarray}

We are interested in computing the steady states, i.e., solutions where $\frac{\partial \mathbf{u}}{\partial t}$ vanishes.
The high-order simulations are obtained through time-advancement, while the ROM solutions are obtained with a fixed-point iteration.

%momentum and pressure equation....

\subsection{Oseen-Iteration}

The \emph{Oseen}-iteration is a secant modulus fixed-point iteration, which in general exhibits a linear rate of convergence \cite{Hess:Oseen}.
Given a current iterate (or initial condition) $\mathbf{u}^k$, the next iterate $\mathbf{u}^{k+1}$ is found
by solving linear system:
\begin{eqnarray}
 -\nu \Delta \mathbf{u}^{k+1} + (\mathbf{u}^k \cdot \nabla) \mathbf{u}^{k+1} + \nabla p &=& \mathbf{f}  \text{ in } \Omega, \label{Hess:eq_Oseen_main} \cl
\nabla \cdot \mathbf{u}^{k+1} &=& 0   \text{ in } \Omega, \cl
 %\mathbf{u} &=& f_d  \text{ on } \partial \Omega,
\mathbf{u}^{k+1} &=& \mathbf{d} \quad \text{ on } \Gamma_D, \cl
\nabla \mathbf{u}^{k+1} \cdot \mathbf{n} &=& \mathbf{g} \quad \text{ on } \Gamma_N. \el
\end{eqnarray}
Iterations are typical stopped when the relative difference between iterates falls below a predefined tolerance in a suitable norm, like the $L^2(\Omega)$ or $H^1_0(\Omega)$ norm. %in the velocity.
%An initial solution $\mathbf{u}^0(\nu_0)$ is computed by time-advancement of \eqref{Hess:NSE0}--\eqref{Hess:NSE1} from zero initial conditions at a parameter value $\nu_0$,
%and the whole parameter domain is then explored by using a continuation method with the \emph{Oseen}-iteration.

%The \emph{Oseen}-iteration is then initialized with $\mathbf{u}(\nu_0)$ to solve for the velocity at a value close to $\nu_0$. 
%Solutions over the whole parameter domain are then explored by using a continuation method, i.e., using 
%the solution of the fixed point iteration for the previous parameter value as initial condition at the next parameter value.
%This improves computation time over initializing the \emph{Oseen}-iteration with a zero field.
  
\subsection{Model Description}

\begin{figure}[t]
\begin{center}
\includegraphics[scale=.18]{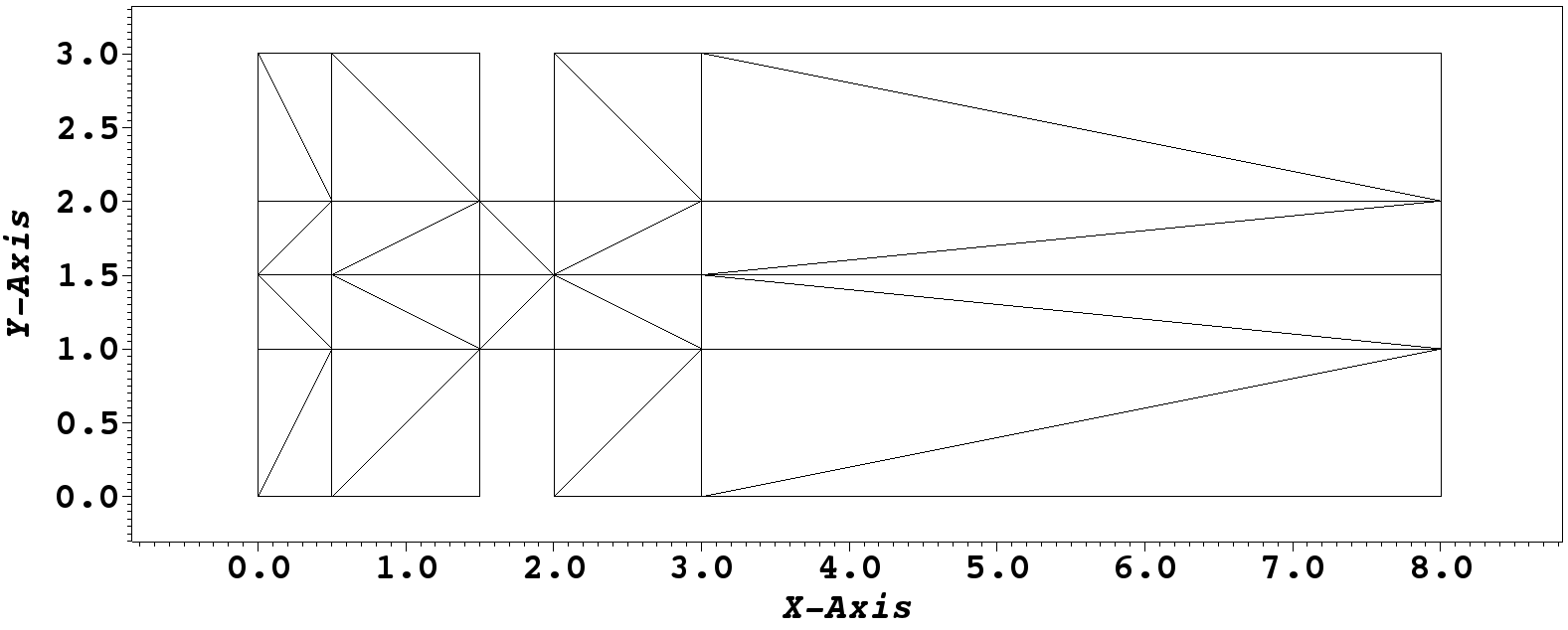}
\caption{Reference computational domain for the channel flow, divided into $36$ triangles. }
\label{Hess:domain}       % Give a unique label
\end{center}
\end{figure}

We consider the reference computational domain shown in Fig.~\ref{Hess:domain}, which is decomposed into $36$ triangular spectral elements. 
The spectral element expansion uses modal Legendre polynomials of the Koornwinder-Dubiner type of order $p = 11$ for the velocity. 
Details on the discretization method can be found in chapter 3.2 of \cite{Hess:Sherwin:2005}.
The pressure \emph{ansatz} space is chosen of order $p-2$ to fulfill the inf-sup stability condition \cite{Hess:BBF,Hess:infsup}.
A parabolic inflow profile is  prescribed at the inlet (i.e., $x = 0$) with horizontal velocity component 
$u_x(0,y) = y(3-y)$ for $y \in [0, 3]$.
At the outlet (i.e $x = 8$) we impose a stress-free boundary condition, everywhere else we prescribe a no-slip condition.

\begin{figure}[ht]
\begin{center}
 \includegraphics[scale=.2]{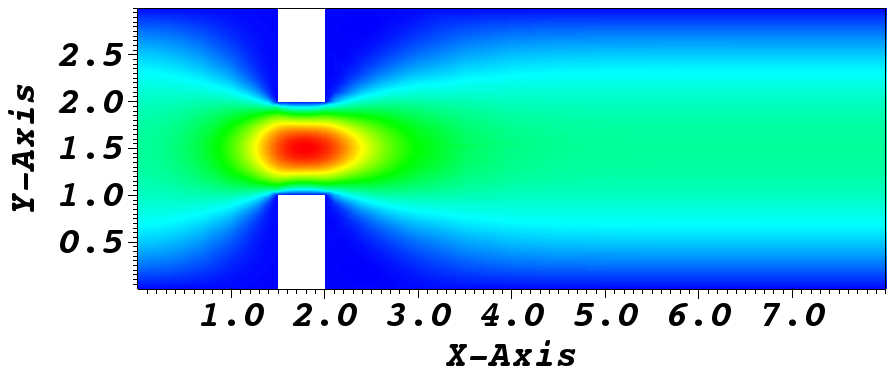} $\quad$
 \includegraphics[scale=.26]{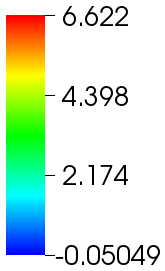} \\
 \includegraphics[scale=.2]{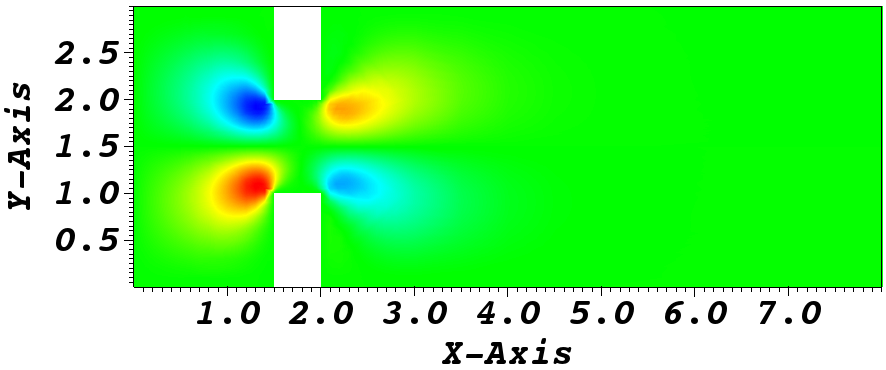} $\quad$
 \includegraphics[scale=.26]{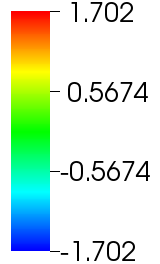}
 \caption{Full order, steady-state solution for $\mu = 1$: velocity in x-direction (top) and y-direction (bottom).}
 \label{Hess:Geo1}
 \end{center}
\end{figure}

\begin{figure}[ht]
\begin{center}
 \includegraphics[scale=.2]{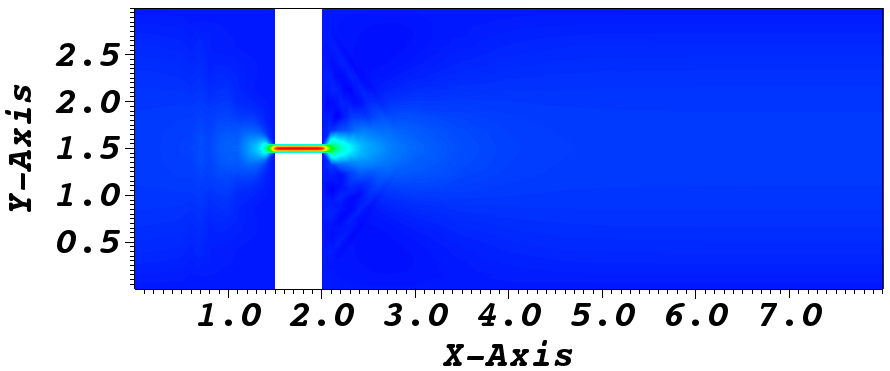} $\quad$
 \includegraphics[scale=.26]{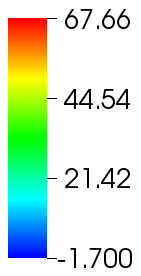} \\
 \includegraphics[scale=.2]{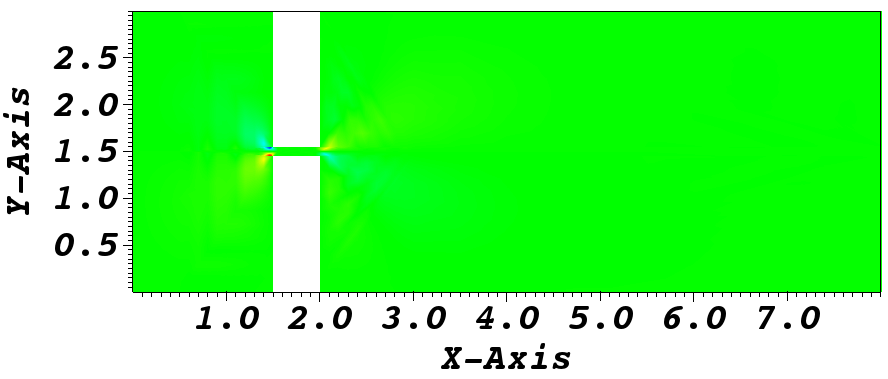} $\quad$
 \includegraphics[scale=.26]{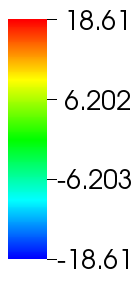}
 \caption{Full order, steady-state solution for $\mu = 0.1$: velocity in x-direction (top) and y-direction (bottom).}
 \label{Hess:Geo0p1}
  \end{center}
\end{figure}

\begin{figure}[ht]
\begin{center}
 \includegraphics[scale=.2]{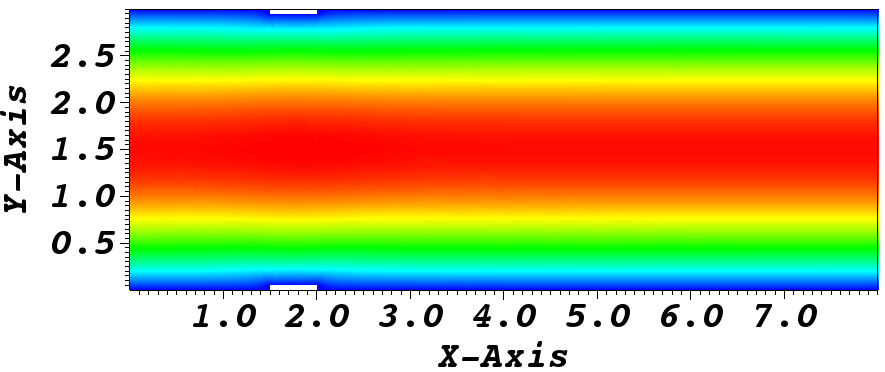} $\quad$
 \includegraphics[scale=.26]{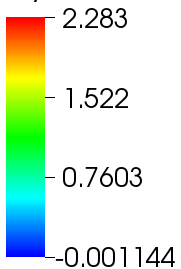} \\
 \includegraphics[scale=.2]{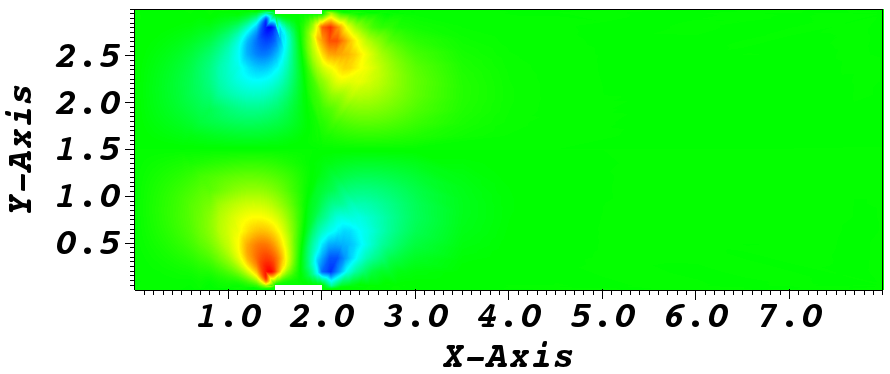} $\quad$
 \includegraphics[scale=.26]{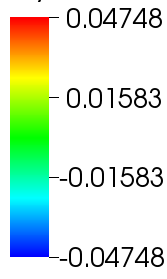}
 \caption{Full order, steady-state solution for $\mu = 2.9$: velocity in x-direction (top) and y-direction (bottom).}
 \label{Hess:Geo2p9}
  \end{center}
\end{figure}

The height of the narrowing in the reference configuration is $\mu = 1$, from $y = 1$ to $y = 2$.
See Fig.~\ref{Hess:domain}.
Parameter $\mu$ is considered variable in the interval $\mu \in [0.1, 2.9]$. 
The narrowing is shrunken or expanded as to maintain
the geometry symmetric about line $y = 1.5$. 
%The high-order field solutions are computed by time advancement and 
Fig.~\ref{Hess:Geo1}, Fig.~\ref{Hess:Geo0p1}, and Fig.~\ref{Hess:Geo2p9}
show the velocity components close to the the steady state for $\mu = 1, 0.1, 2.9$, respectively.


The viscosity is kept constant to $\nu = 1$.
For these simulations, the \emph{Reynolds} number \eqref{eq:re} is between $5$ and $10$,  
with maximum velocity in the narrowing as characteristic velocity $U$
and the height of the narrowing characteristic length $L$. 
For larger \emph{Reynolds} numbers (about $30$), a supercritical pitchfork bifurcation occurs 
giving rise to the so-called \emph{Coanda} effect \cite{Hess:wille_fernholz_1965,Hess:max_bif,Hess:ENUMATH17_me}, which is not subject of the current study.
Our model is similar to the model considered in \cite{Hess:Pitton2017534}, i.e. 
an expansion channel with an inflow profile of varying height. However, in \cite{Hess:Pitton2017534} 
the computational domain itself does not change.

%\begin{eqnarray}
% Re = \frac{1}{4\nu}.
% \label{Hess:Re_visc_example}
%\end{eqnarray}

%Consider a parametric variation in the viscosity $\nu$, ranging from $\nu = 0.0075$ to $\nu = 0.0025$, which corresponds 
%to \emph{Reynolds} numbers between $33$ and $100$.
%The solution for  $\nu = 0.0075$ is shown in Fig.~\ref{Hess:Xp0075}.
%It is slightly unsymmetrical, which marks the onset of the \emph{Coanda} effect \cite{Hess:wille_fernholz_1965}, \cite{Hess:Coanda}, which is a known phenomenon characterized as a `wall-hugging' effect
%ccurring at these \emph{Reynolds} numbers. 
%The solution for  $\nu = 0.0025$ is shown in Fig.~\ref{Hess:Xp0025}.
%Here, the \emph{Coanda} effect is fully developed as the flow orients itself along the boundaries.

\section{Spectral Element Full Order Discretization}

The \emph{Navier-Stokes} problem is discretized with the spectral element 
method. The spectral/hp element software framework used is Nektar++ in version 4.4.0\footnote{See \textbf{www.nektar.info}.}.
%An initial solution is computed with a time-dependent simulation.
%This initial solution then serves as an initial guess for the Oseen-iteration. 
The discretized system of size $N_\delta$ to solve at each step of the \emph{Oseen}-iteration for fixed $\mu$ can be written as
%is given by \eqref{Hess:fully_expanded} as
\begin{eqnarray}
\begin{bmatrix}
\begin{array}{ccc}
A & -D^T_{bnd} & B  \\
-D_{bnd} & 0 & -D_{int} \\
\tilde{B}^T & -D^T_{int} & C
\end{array}
\end{bmatrix}
\begin{bmatrix}
\begin{array}{ccc}
\mathbf{v}_{bnd} \\
\mathbf{p} \\
\mathbf{v}_{int} 
\end{array}
\end{bmatrix}
&=
\begin{bmatrix}
\begin{array}{ccc}
\mathbf{f}_{bnd} \\
\mathbf{0} \\
\mathbf{f}_{int}
\end{array}
\end{bmatrix},
\label{Hess:fully_expanded}
\end{eqnarray}
\noindent where $\mathbf{v}_{bnd}$ and $\mathbf{v}_{int}$ denote velocity degrees of freedom on the boundary 
and in the interior of the domain, respectively, while $\mathbf{p}$
denotes the pressure degrees of freedom. The forcing terms on the boundary and interior
are denoted by $\mathbf{f}_{bnd}$ and $\mathbf{f}_{int}$, respectively.
The matrix $A$ assembles the boundary-boundary coupling, $B$ the boundary-interior coupling, $\tilde{B}$ the interior-boundary coupling, and
$C$ assembles the interior-interior coupling of elemental velocity \emph{ansatz} functions. 
In the case of a \emph{Stokes} system, it holds that $B = \tilde{B}^T$, but this is not the case for the \emph{Oseen} equation
because of the linearized convective term.
%the linearization term $(\mathbf{u}^k \cdot \nabla) \mathbf{u}$ is present in \eqref{Hess:eq_Oseen_main}.
The matrices $D_{bnd}$ and $D_{int}$ assemble the pressure-velocity boundary and pressure-velocity interior contributions, respectively.

The linear system \eqref{Hess:fully_expanded} is assembled in local degrees of freedom, resulting in block matrices $A, B, \tilde{B}, C, D_{bnd}$ and $D_{int}$,
each block corresponding to a spectral element. 
This allows for an efficient matrix assembly since each spectral element is independent from the others, 
but makes the system singular.
In order to solve the system, the local degrees of freedom need to be gathered into the global degrees of freedom \cite{Hess:Sherwin:2005}.

The high-order element solver Nektar++ uses a multilevel static condensation
for the solution of linear systems like \eqref{Hess:fully_expanded}. 
Since static condensation introduces intermediate parameter-dependent matrix inversions (such as $C^{-1}$ in this case) several intermediate projection spaces need to be introduced to use model order reduction \cite{Hess:ENUMATH17_me}.
This can be avoided by instead projecting the expanded system \eqref{Hess:fully_expanded} directly.
The internal degrees of freedom do not need to be gathered, since they are the same in local and global coordinates.
Only \emph{ansatz} functions extending over multiple spectral elements need to be gathered.

%To use the multilevel static condensation in a ROM setting, several intermediate projection spaces need to be introduced \cite{Hess:ENUMATH17_me}.
%This can be avoided by projecting the expanded system directly.
%This is negatively impacting the accuracy as well as very technical from the implementation point of view.
%\anna{Martin: can you please rephrase the previous sentence? I am not sure to understand what it means.}
%We will instead directly project the system \eqref{Hess:fully_expanded}, which is expanded in the local degrees of freedom.

Next, we will take the boundary-boundary coupling across element interfaces into account.
Let $M$ denote the rectangular matrix which gathers the local boundary degrees of freedom into global boundary degrees of freedom.
Multiplication of the first row of \eqref{Hess:fully_expanded} by $M^T M$ will then set the boundary-boundary coupling in local degrees of freedom:
\begin{eqnarray}
\begin{bmatrix}
\begin{array}{ccc}
M^T M A & -M^T M D^T_{bnd} & M^T M B  \\
-D_{bnd} & 0 & -D_{int} \\
\tilde{B}^T & -D^T_{int} & C
\end{array}
\end{bmatrix}
\begin{bmatrix}
\begin{array}{ccc}
\mathbf{v}_{bnd} \\
\mathbf{p} \\
\mathbf{v}_{int} 
\end{array}
\end{bmatrix}
&=
\begin{bmatrix}
\begin{array}{ccc}
M^T M \mathbf{f}_{bnd} \\
\mathbf{0} \\
\mathbf{f}_{int}
\end{array}
\end{bmatrix} .
\label{Hess:fully_expanded_MtM}
\end{eqnarray}

The action of the matrix in \eqref{Hess:fully_expanded_MtM} on the 
degrees of freedom on the Dirichlet boundary is computed and added to the right hand side.
Such degrees of freedom are then removed from \eqref{Hess:fully_expanded_MtM}.
The resulting system can then be used in a projection-based ROM context \cite{Hess:Lassila2014}, 
%in the sense that the projection with snapshot solutions gathers all the local degrees of freedom into a non-singular system matrix 
of high-order dimension $N_\delta \times N_\delta$ and depending on the parameter $\mu$ :
% \begin{eqnarray}
% \begin{bmatrix}
% \begin{array}{ccc}
% \mathcal{A} & -\mathcal{D}^T_{bnd} & \mathcal{B}  \\
% -\mathcal{D}_{bnd} & 0 & -D_{int} \\
% \tilde{\mathcal{B}}^T & -D^T_{int} & C
% \end{array}
% \end{bmatrix}
% \begin{bmatrix}
% \begin{array}{ccc}
% \bar{v}_{bnd} \\
% p \\
% v_{int} 
% \end{array}
% \end{bmatrix}
% &=
% \begin{bmatrix}
% \begin{array}{ccc}
% \bar{f}_{bnd} \\
% 0 \\
% f_{int}
% \end{array}
% \end{bmatrix} .
% \label{Hess:final_to_be_proj}
% \end{eqnarray}
\begin{eqnarray}
\mathcal{A}(\mu) \mathbf{x}(\mu) = \mathbf{f}.
\label{Hess:final_to_be_proj}
\end{eqnarray}

%Note that the velocity boundary degrees of freedom are along the boundaries of the spectral elements and not only the domain boundary, resulting
%in ?? local degrees of freedom for this problem.

%but in a high-order setting this is inferior to the static condensation approach implemented in Nektar++.

%write section 10.1.3 from Nektar user guide

\section{Reduced Order Model}

The reduced order model (ROM) computes accurate approximations to the high-order solutions in the parameter range of interest,
while greatly reducing the overall computational time. This is achieved by two ingredients.
First, a few high-order solutions are computed and the most significant proper orthogonal decomposition (POD) modes are obtained \cite{Hess:Lassila2014}. 
These POD modes define the reduced order ansatz space of dimension $N$, in which the system is solved.
Second, to reduce the computational time, an offline-online computational procedure is used.
See Sec.~\ref{sec:off-on}.
%The affine parameter-dependence allows to precompute all quantities depending on the high-order dimension in an offline phase,
%such that the many-query online phase can be done in negligible time.

The POD computes a singular value decomposition of the snapshot solutions to $99.99\%$ of the most dominant modes \cite{Hess:RBref}, which
define the projection matrix $U \in \mathbb{R}^{N_\delta \times N}$ used
to project system \eqref{Hess:final_to_be_proj}:
\begin{eqnarray}
U^T \mathcal{A}(\mu) U \mathbf{x}_N(\mu) = U^T \mathbf{f} .
\label{Hess:final_proj}
\end{eqnarray}

The low order solution $\mathbf{x}_N(\mu)$ then approximates the high order solution as $\mathbf{x}(\mu) \approx U \mathbf{x}_N(\mu)$.

%The resulting system can then be used in a projection-based ROM context, as is commonly done in 

%such that the reduced order system can be expressed as

%\begin{eqnarray}
%U^T (\hat{A} - \hat{B} \hat{D}^{-1} \hat{C}) U  b_N  = U^T (\hat{f}_{bnd} - \hat{B} \hat{D}^{-1} \hat{f}_p) .
%\label{Hess:hats_very_final_red}
%\end{eqnarray}

%\noindent with the reduced order solution $b_N$.

%While the twice statically condensed system has the lowest dimension, it is difficult to explicitly extract the parameter-dependence from this
%system, which is present in every matrix and required for the offline-online decomposition.
%Instead the offline-online decomposition will proceed from one level of static condensation, i.e., \eqref{Hess:first_condensed}.

\subsection{Offline-Online Decomposition}\label{sec:off-on}

The offline-online decomposition \cite{Hess:RBref} enables the computational speed-up of the ROM approach in many-query scenarios.
It relies on an affine parameter dependency, such that all computations depending on the high-order model size can be moved into
a parameter-independent offline phase, while having a fast input-output evaluation online.

\begin{figure}[t]
\begin{center}
\includegraphics[scale=.23]{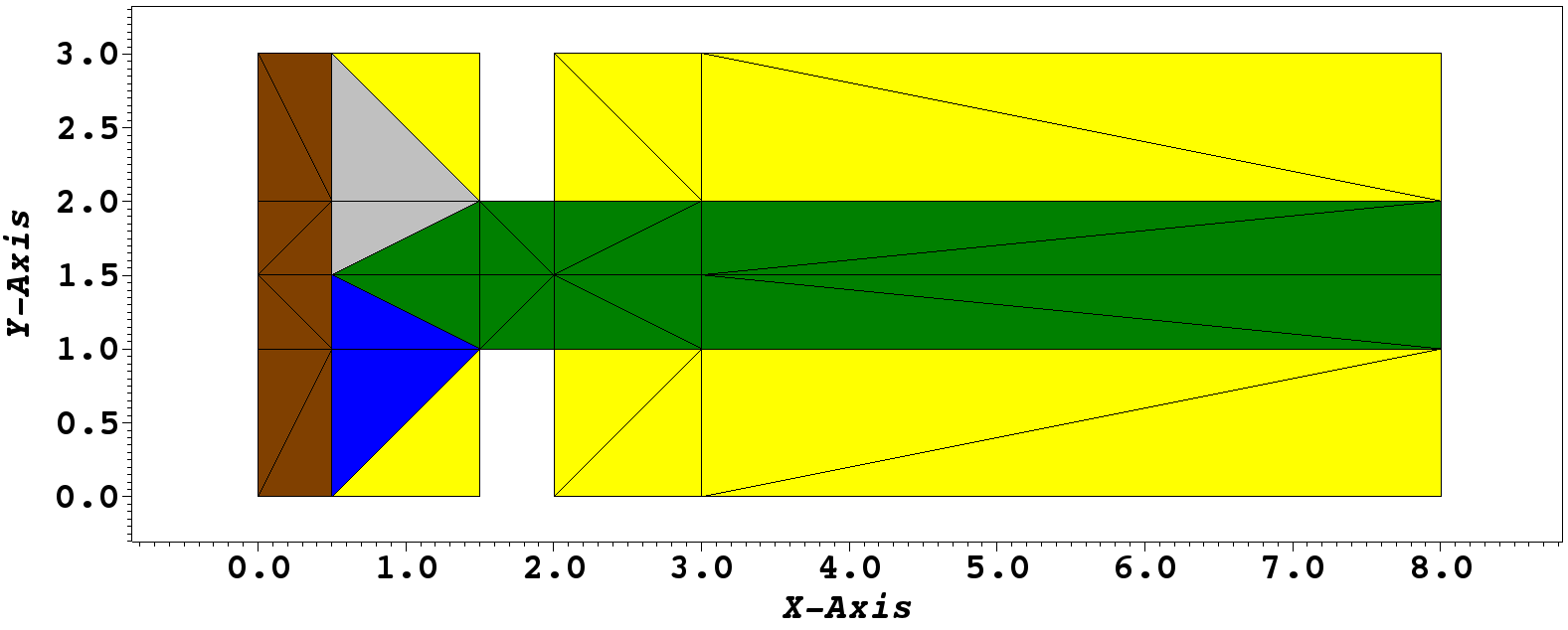}
\caption{Reference computational domain with subdomains $\Omega_1$ (green), $\Omega_2$ (yellow), $\Omega_3$ (blue), $\Omega_4$ (grey) and $\Omega_5$ (brown).}
\label{Hess:subdomain}       % Give a unique label
 \end{center}
\end{figure}

In the example under consideration here, the parameter dependency is already affine and a mapping to the reference domain can be
established without using an approximation technique such as the empirical interpolation method.
Thus, there exists an affine expansion of the system matrix $\mathcal{A}(\mu)$ in the parameter $\mu$ as
\begin{eqnarray}
\mathcal{A}(\mu) = \sum_{i=1}^Q \Theta_i(\mu) \mathcal{A}_i. 
\label{Hess:off-on-gen}
\end{eqnarray}
The coefficients $\Theta_i(\mu)$ are computed from the mapping $\mathbf{x} = T_k(\mu) \hat{\mathbf{x}} + \mathbf{g}_k $, $T_k \in \mathbb{R}^{2 \times 2}$, $\mathbf{g}_k \in \mathbb{R}^{2}$,
 which maps the deformed subdomain $\hat{\Omega}_k$ to
the reference subdomain $\Omega_k$. See also \cite{Hess:Rozza:103010,Hess:Rozza:Q07}. 
Fig.~\ref{Hess:subdomain} shows the reference subdomains $\Omega_k$ for 
the problem under consideration.

For each subdomain $\hat{\Omega}_k$ the elemental basis function evaluations are transformed to the reference domain.
For each velocity basis function $ \mathbf{u} = (u_1, u_2),  \mathbf{v} = (v_1, v_2) ,  \mathbf{w} = (w_1, w_2)  $ and each (scalar) pressure basis function $\psi$, 
we can write the transformation with summation convention as:
\begin{eqnarray}
%\hat{\mathcal{A}} = \sum_{k=1}^K \int_{\hat{\Omega}_k} \frac{\partial \phi_l}{\partial x_i} \hat{\nu}_{ij} \frac{\partial \phi_m}{\partial x_j} d \hat{\Omega}_k \\
%\mathcal{A} = \sum_{k=1}^K \int_{\Omega_k} \frac{\partial \phi_l}{\partial x_i} \nu_{ij} \frac{\partial \phi_m}{\partial x_j} d \hat{\Omega}_k
\int_{\hat{\Omega}_k} \frac{\partial \hat{\mathbf{u}}}{\partial \hat{x}_i} \hat{\nu}_{ij} \frac{\partial \hat{\mathbf{v}}}{\partial \hat{x}_j} d \hat{\Omega}_k 
&=& \int_{\Omega_k} \frac{\partial \mathbf{u}}{\partial x_i} \nu_{ij} \frac{\partial \mathbf{v}}{\partial x_j} d \Omega_k, \cl
\int_{\hat{\Omega}_k} \hat{\psi} \nabla \cdot \hat{\mathbf{u}} d \hat{\Omega}_k 
&=& \int_{\Omega_k} \psi \chi_{ij} \frac{\partial u_j}{\partial x_i} d \Omega_k, \cl
\int_{\hat{\Omega}_k} (\hat{\mathbf{u}} \cdot \nabla) \hat{\mathbf{v}} \cdot \hat{\mathbf{w}} d \hat{\Omega}_k 
&=& \int_{\Omega_k} u_i \pi_{ij} \frac{\partial v_j}{\partial x_i} \mathbf{w} d \Omega_k, \el
\label{Hess:trafo_terms}
\end{eqnarray}
\noindent with 
\begin{eqnarray}
 \nu_{ij} &=& T_{ii'} \hat{\nu}_{i'j'} T_{jj'} \det(T)^{-1},  \cl
 \chi_{ij} &=& \pi_{ij} =  T_{ij} \det(T)^{-1}. \el
\end{eqnarray}

The subdomain $\Omega_5$ (see Fig.~\ref{Hess:subdomain})
is kept constant, so that no interpolation of the inflow profile is necessary.
To achieve fast reduced order solves, the offline-online decomposition expands 
the system matrix as in \eqref{Hess:off-on-gen} and
%in the parameter of interest and 
computes the parameter independent projections offline, which are
stored as small-sized matrices of the order $N \times N$.
Since in an \emph{Oseen}-iteration each matrix is dependent on the previous iterate, the submatrices corresponding to each basis function are assembled and then 
formed online using the reduced basis coordinate representation of the current iterate. 
This is the same procedure used for the assembly of the nonlinear term in the \emph{Navier-Stokes} case \cite{Hess:Lassila2014}.

%This is made explicit here on the example of the matrix $A$ and analogous for the other cases.
%Since $A$ denotes the boundary-boundary coupling of velocity degrees of freedom, each element $(i,j)$ of $A$ is of the form

%\begin{eqnarray}
% A_{ij} = (\nabla \Phi^b_i, \nu \nabla \Phi^b_i) + (\nabla \Phi^b_i, u_k \cdot \nabla \Phi^b_j),
%\end{eqnarray}

%\noindent with spectral element \emph{ansatz} function on the boundary $\Phi^b_i$ and current \emph{Oseen} iterate $u_k$.
%Since the viscosity is considered a parameter, $\nu$ is separated and the parameter-independent part $(\nabla \Phi^b_i, \nabla \Phi^b_i)$ is 
%assembled offline.

%introduce reduced quantities here

%snapshot space generation

%\subsection{Reduced order Model within static condensation}

%A further gain can be achieved by not only using the ROM to speed up the computation time for the linear system \eqref{Hess:hats_final}, but
%also when inverting the block matrices in $D$ and $\hat{D}$.

%A direct approach did not yield good results...

%instead lets do this special thingy here....

\section{Numerical Results}

\begin{figure}[ht]
 \center
 \includegraphics[scale=.7]{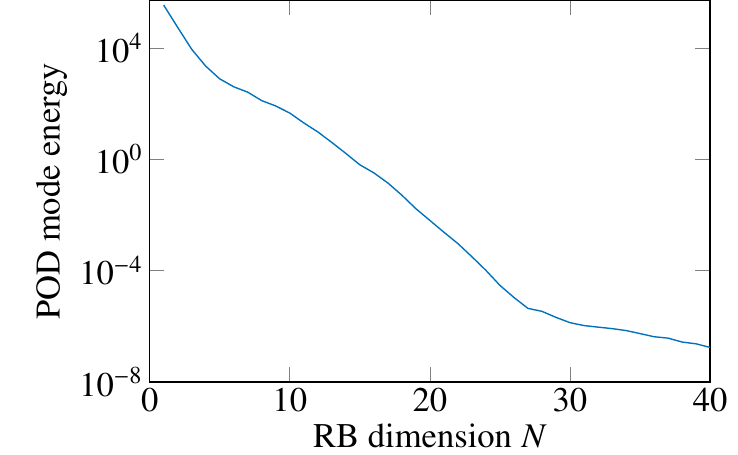}
% \caption{Decay of POD mode energy.}
 \includegraphics[scale=.7]{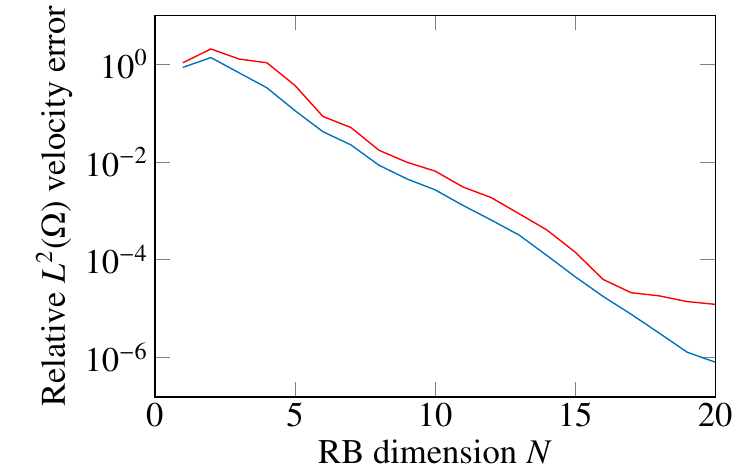}
 \caption{Left: Decay of POD mode energy. Right: Maximum (red) and mean (blue) relative $L^2(\Omega)$ error for the velocity over increasing reduced basis dimension.}
 \label{Hess:relerr}
\end{figure}

The accuracy of the ROM is assessed using $40$ snapshots sampled 
uniformly over the parameter domain $[0.1, 2.9]$ for the POD and 
$40$ randomly chosen parameter locations to test the accuracy.
Fig.~\ref{Hess:relerr} (left) shows the decay of the energy of the POD modes. 
To reach the typical threshold of $99.99\%$ on the POD energy, it
takes $9$ POD modes as RB ansatz functions.
Fig.~\ref{Hess:relerr} (right) shows the relative $L^2(\Omega)$ approximation error of the reduced order model with respect to the full order model
up to $6$ digits of accuracy, evaluated at the $40$ randomly chosen verification parameter locations. 
With $9$ POD modes the maximum approximation error is less than $0.7\%$ 
and the mean approximation error is less than $0.5\%$.

While the full-order solves were computed with Nektar++, the reduced-order computations were done in ITHACA-SEM 
with a separate python code.
To assess the computational gain, the time for
a fixed point iteration step using the full-order system is compared to the time for
a fixed point iteration step of the ROM with dimension $20$, both done in python.
The ROM online phase reduces the computational time by a factor of over $100$.
The offline time is dominated by computing the snapshots and 
the trilinear forms used to project the advection terms. See \cite{Hess:Lassila2014} for detailed explanations.

%model description with picture and one snapshot -- maybe later

%do the comparison

%expect two graphs 
%- fast approx with standard large RB
%- faster approx with 'statically condensed RB'
%  maybe additional cost of element wise C-solves

\section{Conclusion and Outlook}

We showed that the POD reduced basis technique generates accurate reduced order models for SEM discretized models under parametric variation of the geometry.
The potential of a high-order spectral element method 
with a reduced basis ROM is the subject of current investigations. See also \cite{Hess:Taddei}.
Since each spectral element comprises a block in the system matrix in local coordinates, 
a variant of the reduced basis element 
method (RBEM) (\cite{Hess:Maday2002195,Hess:lovgren_maday_ronquist_2006}) 
can be successfully applied in the future.

%when possible to statically condensed RB
%LINK with RBEM ??

\section*{Acknowledgments}  

This work was supported by European Union Funding for Research and Innovation through the European Research Council
(project H2020 ERC CoG 2015 AROMA-CFD project 681447, P.I. Prof. G. Rozza).
This work was also partially supported by NSF through grant DMS-1620384 (Prof. A. Quaini).


\begin{thebibliography}{99.}%
% Use the following syntax and markup for your references if 
% the subject of your book is from the field 
% "Mathematics, Physics, Statistics, Computer Science"
%
% Contribution 
%\bibitem{science-contrib} Broy, M.: Software engineering --- from auxiliary to key technologies. In: Broy, M., Dener, E. (eds.) Software Pioneers, pp. 10-13. Springer, Heidelberg (2002)
%
% Online Document
%\bibitem{science-online} Dod, J.: Effective substances. In: The Dictionary of Substances and Their Effects. Royal Society of Chemistry (1999) Available via DIALOG. \\
%\url{http://www.rsc.org/dose/title of subordinate document. Cited 15 Jan 1999}
%
% Monograph
%\bibitem{science-mono} Geddes, K.O., Czapor, S.R., Labahn, G.: Algorithms for Computer Algebra. Kluwer, Boston (1992) 
%
% Journal article
%\bibitem{science-journal} Hamburger, C.: Quasimonotonicity, regularity and duality for nonlinear systems of partial differential equations. Ann. Mat. Pura. Appl. \textbf{169}, 321--354 (1995)
%
% Journal article by DOI
%\bibitem{science-DOI} Slifka, M.K., Whitton, J.L.: Clinical implications of dysregulated cytokine production. J. Mol. Med. (2000) doi: 10.1007/s001090000086 
%


\bibitem{Hess:Sherwin:2005} 
Karniadakis, G., Sherwin, S.:
Spectral/hp Element Methods for Computational Fluid Dynamics.
In: Oxford University Press, 2nd ed. (2005).


\bibitem{Hess:CHQZ1}
Canuto, C., Hussaini, M.Y., Quarteroni, A., Zhang, Th.A.:
Spectral {M}ethods {F}undamentals in {S}ingle {D}omains.
In: Springer -- Scientific {C}omputation, (2006).


\bibitem{Hess:CHQZ2}
Canuto, C., Hussaini, M.Y., Quarteroni, A., Zhang, Th.A.:
Spectral {M}ethods {E}volution to {C}omplex {G}eometries and {A}pplications to {F}luid {D}ynamics.
In: Springer -- Scientific {C}omputation, (2007).


\bibitem{Hess:PateraSEM}
Patera, A.T.:
A Spectral Element Method for Fluid Dynamics; Laminar Flow in a Channel Expansion.
Journal of Computational Physics, {\bf 54}:3 (1984), 468--488.

\bibitem{Hess:Herrero2013132}
Herrero, H., Maday, Y., Pla, F.:
RB (Reduced Basis) for RB (Rayleigh--B\'enard).
Computer Methods in Applied Mechanics and Engineering, {\bf 261--262}, (2013), 132--141.

\bibitem{Hess:Taddei}
Fick, L., Maday, Y., Patera A., Taddei T.:
A stabilized POD model for turbulent flows over a range of Reynolds numbers: Optimal parameter sampling and constrained projection.
Journal of Computational Physics, {\bf 371} (2018), 214--243.


\bibitem{Hess:nektar}
 Cantwell, C.D., Moxey, D. , Comerford, A., Bolis, A., Rocco, G., Mengaldo, G., de Grazia, D., Yakovlev, S., Lombard, J.-E., Ekelschot, D., Jordi, B., Xu, H., Mohamied, Y., Eskilsson, C., Nelson, B., Vos, P., Biotto, C., Kirby, R.M., Sherwin, S.J.:
 Nektar++: An open-source spectral/hp element framework.
Computer Physics Communications, {\bf 192}, (2015), 205--219.


\bibitem{Hess:Holmes}
Holmes, P., Lumley, J., Berkooz, G.:
Turbulence, Coherent Structures, Dynamical Systems and Symmetry.
Cambridge University Press,  (1996).

\bibitem{Hess:Pitton2017534}
Pitton, G., Quaini, A., Rozza, G.:
Computational Reduction Strategies for the Detection of Steady Bifurcations in Incompressible Fluid-Dynamics: Applications to \emph{Coanda} Effect in Cardiology.
Journal of Computational Physics, {\bf 344}, (2017), 534--557.




\bibitem{Hess:RBref}
Hesthaven, J.S., Rozza, G., Stamm, B.:
Certified Reduced Basis Methods for Parametrized Partial Differential Equations.
In: SpringerBriefs in Mathematics, (2016).

\bibitem{Hess:Oseen}
Burger, M.:
Numerical Methods for Incompressible Flow, Lecture Notes.
UCLA (2010).




%\bibitem{Hess:dijkstra_wubs_cliffe_doedel_dragomirescu_eckhardt_gelfgat_hazel_lucarini_salinger_etal._2014}
%{\sc Dijkstra, Henk A. and Wubs, Fred W. and Cliffe, Andrew K. and Doedel, Eusebius and Dragomirescu, Ioana F. and Eckhardt, Bruno and Gelfgat, Alexander Yu. and Hazel, Andrew L. and Lucarini, Valerio and Salinger, Andy G. and et al.}, 
%{\em Numerical Bifurcation Methods and their Application to Fluid Dynamics: Analysis beyond Simulation}, 
%Communications in Computational Physics, {\bf 15}: 1, 
%(2014), 1--45.






%\bibitem{Hess:podcvtlang}
%{\sc S. Ullmann and J. Lang},
%{\em POD and CVT Galerkin reduced-order modelling of the flow around a cylinder},
%PAMM, (2012), {\bf 12}.

\bibitem{Hess:BBF}
Boffi, D., Brezzi F., Fortin, M.:
Mixed Finite Element Methods and Applications.
Springer Series in Computational Mathematics, (2013).


\bibitem{Hess:infsup}
Quarteroni, A., Valli, A.:
Numerical  Approximation  of  Partial  Differential  Equations.
Springer-Verlag, Berlin-Heidelberg, (1994).




\bibitem{Hess:Pitton2017}
Pitton, G., Rozza, G.:
On the Application of Reduced Basis Methods to Bifurcation Problems in Incompressible Fluid Dynamics.
Journal of Scientific Computing, (2017).

\bibitem{Hess:wille_fernholz_1965}
Wille, R., Fernholz, H.: 
Report on the first European Mechanics Colloquium, on the \emph{Coanda} effect.
Journal of Fluid Mechanics, {\bf 23}:4 (1965), 801--819.

%\bibitem{Hess:Coanda}
%Quaini, A., Glowinski, R., {\v Cani\'c}, S.:
%A  computational  study  on  the  generation  of  the  \emph{Coanda} effect in a mock heart chamber.
%RIMS K\^oky\^uroku series, No. 2009-4, (2016).






%\bibitem{Hess:10.1080/10618562.2016.1144877}
%{\sc A. Quaini, R. Glowinski and S. Čanić},
%{\em Symmetry Breaking and Preliminary Results about a Hopf Bifurcation for Incompressible Viscous Flow in an Expansion Channel},
%International Journal of Computational Fluid Dynamics, {\bf 30}:1
%(2016), 7--19.

\bibitem{Hess:Rozza:103010}
Rozza, G.:
Real-time reduced basis solutions for Navier-Stokes  equations: optimization of parametrized bypass  configurations.
ECCOMAS CFD 2006 Proceedings on CD, 676 (2006), 1--16.

\bibitem{Hess:Rozza:Q07}
Quarteroni, A., Rozza, G.:
Numerical   Solution   of   parametrized   Navier-Stokes equations by reduced basis methods. 
Num. Meth. Part. Diff. Eq., 23(4) (2007), 923--948.


\bibitem{Hess:Lassila2014}
Lassila, T., Manzoni, A., Quarteroni, A., Rozza, G.:
Model Order Reduction in Fluid Dynamics: Challenges and Perspectives.
In: Reduced Order Methods for Modelling and Computational Reduction, Springer International Publishing, MS\&A, Vol. 9, A. Quarteroni, G.Rozza eds. (2014), 235--273.



%\bibitem{Hess:Yano20130036}
%{\sc M.~Yano and A.~T.~Patera}, 
%{\em A Space{\textendash}Time Variational Approach to Hydrodynamic Stability Theory},
%Proceedings of the Royal Society of London A: Mathematical, Physical and Engineering Sciences {\bf 469}:2155
% (2013).


\bibitem{Hess:Maday2002195}
Maday, Y., Ronquist, E.M.: 
A Reduced-Basis Element Method.
Comptes Rendus Mathematique {\bf 335}:2 (2002), 195--200.


\bibitem{Hess:lovgren_maday_ronquist_2006}
Lovgren, A.E., Maday, Y, Ronquist, E.M.:
A Reduced Basis Element Method for the Steady Stokes Problem.
ESAIM: Mathematical Modelling and Numerical Analysis {\bf 40}:3  (2006), 529--552.



\bibitem{Hess:ENUMATH17_me}
Hess, M.W., Rozza, G.:
A Spectral Element Reduced Basis Method in Parametric {CFD}.
Numerical Mathematics and Advanced Applications - ENUMATH 2017, Springer, in press, ArXiv e-print 1712.06432, (2018).

\bibitem{Hess:max_bif}
Hess, M.W., Alla, A., Quaini, A., Rozza, G., Gunzburger, M.: 
A Localized Reduced-Order Modeling Approach for PDEs with Bifurcating Solutions.
ArXiv e-print 1807.08851, accepted for publication in Computer Methods in Applied Mechanics and Engineering (CMAME), (2019).




\end{thebibliography}
\end{document}